\documentclass[a4paper,11pt]{article}
\usepackage{fullpage}
\usepackage{graphics}
\usepackage[latin2]{inputenc}
\usepackage{multicol}
\usepackage{amsmath,amssymb,amsbsy,amsthm,amsfonts}
\usepackage{epsfig} 

\begin{document}
\def\note#1{\marginpar{\small #1}}

\def\tens#1{\pmb{\mathsf{#1}}}
\def\vec#1{\boldsymbol{#1}}

\def\norm#1{\left|\!\left| #1 \right|\!\right|}
\def\fnorm#1{|\!| #1 |\!|}
\def\abs#1{\left| #1 \right|}
\def\ti{\text{I}}
\def\tii{\text{I\!I}}
\def\tiii{\text{I\!I\!I}}

\def\diver{\mathop{\mathrm{div}}\nolimits}
\def\grad{\mathop{\mathrm{grad}}\nolimits}
\def\Div{\mathop{\mathrm{Div}}\nolimits}
\def\Grad{\mathop{\mathrm{Grad}}\nolimits}

\def\tr{\mathop{\mathrm{tr}}\nolimits}
\def\cof{\mathop{\mathrm{cof}}\nolimits}
\def\det{\mathop{\mathrm{det}}\nolimits}

\def\lin{\mathop{\mathrm{span}}\nolimits}
\def\pr{\noindent \textbf{Proof: }}
\def\pp#1#2{\frac{\partial #1}{\partial #2}}
\def\dd#1#2{\frac{\d #1}{\d #2}}

\def\T{\mathcal{T}}
\def\R{\mathbb{R}}
\def\bx{\vec{x}}
\def\be{\vec{e}}
\def\bef{\vec{f}}
\def\bec{\vec{c}}
\def\bs{\vec{s}}
\def\ba{\vec{a}}
\def\bn{\vec{n}}
\def\bphi{\vec{\varphi}}
\def\btau{\vec{\tau}}
\def\bc{\vec{c}}
\def\bg{\vec{g}}

\def\bW{\vec{W}}
\def\bT{\tens{T}}
\def\bD{\tens{D}}
\def\bF{\tens{F}}
\def\bB{\tens{B}}
\def\bV{\tens{V}}
\def\bS{\tens{S}}
\def\bI{\tens{I}}
\def\bi{\vec{i}}
\def\bv{\vec{v}}
\def\bfi{\vec{\varphi}}
\def\bk{\vec{k}}
\def\b0{\vec{0}}
\def\bom{\vec{\omega}}
\def\bw{\vec{w}}
\def\p{\pi}
\def\bu{\vec{u}}

\def\ID{\mathcal{I}_{\bD}}
\def\IP{\mathcal{I}_{p}}
\def\Pn{(\mathcal{P})}
\def\Pe{(\mathcal{P}^{\eta})}
\def\Pee{(\mathcal{P}^{\varepsilon, \eta})}

\def\Ln#1{L^{#1}_{\bn}}

\def\Wn#1{W^{1,#1}_{\bn}}

\def\Lnd#1{L^{#1}_{\bn, \diver}}

\def\Wnd#1{W^{1,#1}_{\bn, \diver}}

\def\Wndm#1{W^{-1,#1}_{\bn, \diver}}

\def\Wnm#1{W^{-1,#1}_{\bn}}

\def\Lb#1{L^{#1}(\partial \Omega)}

\def\Lnt#1{L^{#1}_{\bn, \btau}}

\def\Wnt#1{W^{1,#1}_{\bn, \btau}}

\def\Lnd#1{L^{#1}_{\bn, \btau, \diver}}

\def\Wntd#1{W^{1,#1}_{\bn, \btau, \diver}}

\def\Wntdm#1{W^{-1,#1}_{\bn,\btau, \diver}}

\def\Wntm#1{W^{-1,#1}_{\bn, \btau}}

\newtheorem{Theorem}{Theorem}[section]
\newtheorem{Example}{Example}[section]
\newtheorem{Lem}{Lemma}[section]
\newtheorem{Rem}{Remark}[section]
\newtheorem{Def}{Definition}[section]
\newtheorem{Col}{Corollary}[section]
\newtheorem{Proposition}{Proposition}[section]

\newcommand{\Om}{\Omega}
\newcommand{ \vit}{\hbox{\bf u}}
\newcommand{ \Vit}{\hbox{\bf U}}
\newcommand{ \vitm}{\hbox{\bf w}}
\newcommand{ \ra}{\hbox{\bf r}}
\newcommand{ \vittest }{\hbox{\bf v}}
\newcommand{ \wit}{\hbox{\bf w}}
\newcommand{ \fin}{\hfill $\square$}

\newcommand{\ZZ}{\mathbb{Z}}
\newcommand{\CC}{\mathbb{C}}
\newcommand{\NN}{\mathbb{N}}
\newcommand{\V}{\zeta}
\newcommand{\RR}{\mathbb{R}}
\newcommand{\EE}{\varepsilon}
\newcommand{\Lip}{\textnormal{Lip}}
\newcommand{\XX}{X_{t,|\textnormal{D}|}}
\newcommand{\PP}{\mathfrak{p}}
\newcommand{\VV}{\bar{v}_{\nu}}
\newcommand{\QQ}{\mathbb{Q}}
\newcommand{\HH}{\ell}
\newcommand{\MM}{\mathfrak{m}}
\newcommand{\rr}{\mathcal{R}}
\newcommand{\tore}{\mathbb{T}_3}
\newcommand{\Z}{\mathbb{Z}}
\newcommand{\N}{\mathbb{N}}

\newcommand{\F}{\overline{\boldsymbol{\tau} }}

\newcommand{\moy} {\overline {\vit} }
\newcommand{\moys} {\overline {u} }
\newcommand{\mmoy} {\overline {\wit} }
\newcommand{\g} {\nabla }
\newcommand{\G} {\Gamma }
\newcommand{\x} {{\bf x}}
\newcommand{\E} {\varepsilon}
\newcommand{\BEQ} {\begin{equation} }
\newcommand{\EEQ} {\end{equation} }
\makeatletter
\@addtoreset{equation}{section}
\renewcommand{\theequation}{\arabic{section}.\arabic{equation}}

\newcommand{\hs}{{\rm I} \! {\rm H}_s}
\newcommand{\esp} [1] { {\bf I} \! {\bf H}_{#1} }

\newcommand{\vect}[1] { \overrightarrow #1}

\newcommand{\hsd}{{\rm I} \! {\rm H}_{s+2}}

\newcommand{\HS}{{\bf I} \! {\bf H}_s}
\newcommand{\HSD}{{\bf I} \! {\bf H}_{s+2}}

\newcommand{\hh}{{\rm I} \! {\rm H}}
\newcommand{\lp}{{\rm I} \! {\rm L}_p}
\newcommand{\leb}{{\rm I} \! {\rm L}}
\newcommand{\lprime}{{\rm I} \! {\rm L}_{p'}}
\newcommand{\ldeux}{{\rm I} \! {\rm L}_2}
\newcommand{\lun}{{\rm I} \! {\rm L}_1}
\newcommand{\linf}{{\rm I} \! {\rm L}_\infty}
\newcommand{\expk}{e^{ {\rm i} \, { k_3} x_3}}
\newcommand{\proj}{{\rm I}Ê\! {\rm P}}

\renewcommand{\theenumi}{\Roman{section}.\arabic{enumi}}

\newcounter{taskcounter}[section]

\newcommand{\bib}[1]{\refstepcounter{taskcounter} {\begin{tabular}{ l p{13,5cm}} \hskip -0,2cm [\Roman{section}.\roman{taskcounter}] & {#1}
\end{tabular}}}

\renewcommand{\thetaskcounter}{\Roman{section}.\roman{taskcounter}}

\newcounter{technique}[section]

\renewcommand{\thetechnique}{\roman{section}.\roman{technique}}

\newcommand{\tech}[1]{\refstepcounter{technique} {({\roman{section}.\roman {technique}}) {\rm  #1}}}

\newcommand{\B}{\mathcal{B}}

\newcommand{\diameter}{\operatorname{diameter}}


%
%
%
%


%
%

\begin{center}
\Large{\textbf{ Approximate Deconvolution Model  in a
bounded  domain with a vertical regularization}}
\end{center}
\begin{center}
\textbf{{Hani Ali}}$^{\hbox{a}}$\\
$^{\hbox{a}}$\textit{MAP5, CNRS UMR 8145,\\
 Université Paris Descartes,\\
 75006 Paris,
France\\
hani.ali@parisdescartes.fr}
\end{center}
\bigskip
\textbf{Abstract}
\bigskip\\
We study the global existence issue for a three-dimensional Approximate Deconvolution Model  with a
vertical filter. We  consider this model in a bounded cylindrical domain where  we construct a unique global weak solution. The proof is based on a refinement of the energy method given by Berselli in \cite{Berselli2012}. 

\smallskip
\smallskip

\hskip-0.45cm  \textbf{ MSC: 76D05; 35Q30; 76F65; 76D03}\\


\section{{Introduction and notation}}

The Large eddy simulation model (LES) with anisotropic regularization is introduced in \cite{Berselli2012} in order to study flows in  bounded domains. \\
Instead of choosing  the classical filter \begin{align}
\label{exact}
\mathbb{A} := I- \left( \alpha_1^2 \partial_1^2 + \alpha_2^2 \partial_2^2+  \alpha_3^2  \partial_3^2 \right), \end{align}
the author in \cite{Berselli2012}, considers a horizontal filter
 \begin{align}
\label{horizontal}
\mathbb{A}_h := I-\left( \alpha_1^2 \partial_1^2 + \alpha_2^2 \partial_2^2 \right). \end{align}
This  filter 
is less memory consuming than the classical one \cite{germano, FDT02, LL03}. Moreover, there is no need to introduce
artificial boundary conditions for  Helmholtz operator. However,  
 the smoothing created by  both of these filters  is unnecessary strong. For this reason one can replace  the classical Laplace operator  by  the Laplace operator with fractional regularization $\theta$ and seek for the limiting case where we can
prove global existence and uniqueness of regular solutions.  
In a series of papers \cite{OT2007,bresslilewandowski,A01,HA11b}, it was shown that the LES models which are derived by using instead of the operator (\ref{exact}), the operator
$\mathbb{A}_{\theta}= I + \alpha^{2\theta}(-\Delta)^{\theta}$ for suitable (large enough) $\theta< 1$, are still well-posed.

In this paper, we  consider  LES models with fractional  filter acting only in  one variable
 \begin{align}
\label{vertical}
\mathbb{A}_{3,\theta} := I + \alpha_3^{2\theta} {(-\partial_{3})^{2\theta}}, \quad 0 \le  \theta \le 1. \end{align}
In particular,  we study the global existence and uniqueness of
solutions to the vertical LES model  on a bounded product
domain of the type $D=\Omega \times (-\pi,\pi)$, where $ \Omega $ 
 is a smooth domain, with homogeneous Dirichlet boundary conditions on the lateral boundary $\partial \Omega \times (-\pi,\pi)$, and with periodic boundary conditions in the vertical variable.
For simplicity, we consider the domain  $D=\left\{ \vec{x} \in \R^3, x_1^2 + x_2^2 < d,  -\pi < x_3 < \pi  \right\}$ with $2\pi$ periodicity with respect to $ x_3$.
 More precisely,
we consider the system of equations 

\begin{equation}
\label{alpha ns}
 \left\{
 \begin{array} {llll} \displaystyle
 \partial_t\bw + \overline{\nabla \cdot \bw  \otimes \bw } -  \nu \Delta \bw + \nabla q
 =  \overline{\bef}, \textrm{ in } D, \  t >0,\\
 \nabla \cdot  \bw_{} =0, \textrm{ in } D, \ t >0,\\
 \bw\vert_{\partial \Omega \times (- \pi , \pi)} = 0,  \ t >0,  \\
\bw(x_1,x_2,x_3)= \bw(x_1,x_2,x_3 + 2\pi), \  t >0,\\
\bw(\bx, 0) = \bw_{0}(\bx) =\overline{\bv_{0}}.
 \end{array}
\right.
\end{equation}
 We recall that, in the above equations, the symbol $``\displaystyle{}^{\overline{\quad}}"$  
denotes the vertical filter  (\ref{vertical}), applied component-by-component to the various tensor fields.
 Given that the filter is acting only in the  vertical variable, it is  possible
to require the periodicity only in $x_3$. Moreover, we consider  the filtered function with  homogeneous Dirichlet boundary conditions on the  lateral boundary $\partial \Omega \times (-\pi,\pi)$. These boundary conditions  of the filtered function are supposed to be the same as the unfiltered ones, in order to  prevent from introducing  artificial
boundary conditions.

The use of the vertical filter can be extended to a family of 
 Approximate Deconvolution model.   The   deconvolution family including (\ref{alpha ns}) as the zeroth order case, is
given by

\begin{equation}
\label{dalpha ns}
 \left\{
 \begin{array} {llll} \displaystyle
 \partial_t\bw + \overline{\nabla \cdot D_{N,\theta}(\bw)  \otimes D_{N,\theta}(\bw) } -  \nu \Delta \bw + \nabla q
 =  \overline{\bef}, \textrm{ in } D, \  t >0,\\
 \nabla \cdot  \bw_{} =0, \textrm{ in } D, \ t >0,\\
 \bw\vert_{\partial \Omega \times (- \pi , \pi)} = 0,  \ t >0,  \\
\bw(x_1,x_2,x_3)= \bw(x_1,x_2,x_3 + 2\pi), \  t >0 ,\\
\bw(\bx, 0) = \bw_{0}(\bx) =\overline{\bv_{0}},
 \end{array}
\right.
\end{equation}
where $D_{N,\theta}$  is a deconvolution operator of order N and is given by 
 \begin{align}
D_{N,\theta} = \sum_{i = 0}^N ( I - \mathbb{A}_{3,\theta}^{-1} )^i, \quad  0 \le \theta \le 1.
\end{align}


In order to state our main result pertaining to System (\ref{dalpha ns}), let us introduce our notation.\\
We set:
\begin{align}\vec{L}^2(D):=\left\{ \vec{v}  \in L^2( D)^3, 2\pi\textrm{-periodic in }x_3, \textrm{and }\int_{-\pi}^{\pi} \vec{v} dx_3 = 0   \right\},\end{align}
endowed with the norm
$\| \vec{v}\|_{2}:= \displaystyle \int_{D} \vec{v} \cdot \vec{v} \ d \vec{x} $ which is the usual norm in $L^2(D)^3$.\\
Then, we define the following spaces that are commonly used in the study of the NSE: 
\begin{align}
\vec{H}:= \left\{  \vec{v} \in   \vec{L}^2(D),\textrm{ such that } \nabla \cdot \vec{v} = 0 \textrm{ and }\vec{v} \cdot n = 0 \textrm{ on  }\partial D    \right\},\\
\vec{V}:= \left\{  \vec{v} \in   \vec{H}, \textrm{ such that } \nabla \vec{v} \in  \vec{L}^2(D) \textrm{ and }\vec{v}  = 0 \textrm{ on  }\partial D    \right\}.
\end{align}
Our main result is the following.

 \begin{Theorem} 
Assume  ${\bef} \in L^{2}(0,T;\vec{H}^{}_{})$ be a divergence free function and  $\bv_0 \in \vec{H}$.  let $ 0 \le N <\infty  $ be a given and fixed number and let $\theta > \frac{1}{2}$.  Then problem  (\ref{dalpha ns}), including   problem (\ref{alpha ns}) when $N= 0$, has a unique regular
weak solution.
\label{TH1}
\end{Theorem}

This result is one of the few results that consider the LES model  on a bounded  domain.  It holds also true on the whole space $\R^3$ and  on the torus $\tore$.   Other $\alpha$ models, with partial filter,  will be reported in a forthcoming paper.

The remaining part of the paper is organized as follows. In the second section,  we introduce
several notations and spaces. We also prove auxiliary results which will be
used in the third section.  In the third section, we  prove Theorem \ref{TH1}. Finally, a functional inequality is
given  in the appendix.

\section{ Preliminaries and auxiliary results}

Let $1 < p \le + \infty $ and $1 < q \le + \infty.$ We denote by $\displaystyle L^q_v
L^p_h(D)=
L^q\left((-\pi, +\pi);L^p(\Omega
)\right)$ the space of functions $g$ such
that 
$\displaystyle\int_{-\pi}^{+\pi}\left(\int_{\Omega}|g(x_1, x_2, x_3)|^{p}dx_1dx_2)^{q/p} dx_3\right)^{
1/q} < +\infty. $
In what follows, we denote by $\nabla_h$ the gradient operator with respect to the variables $x_1$ and $x_2$. Of course,  we have 
$\| \nabla_h \vec{v}\|_{2} \le \| \nabla \vec{v}\|_{2} $, for any $ \vec{v} \in \vec{V}$.\\ 
Next, we will give some preliminary results  which will play a important role in the proof of the main theorem.

\begin{Lem}
\label{paicuraugel}
For any $\vec{u} \in L^{2}(D)^3$ satisfying homogeneous
Dirichlet boundary conditions on the boundary $\partial \Omega \times (-\pi,\pi)$ , such that $ \nabla_h \vec{u} \in L^{2}(D)^3$, we have the following 
estimate

\begin{align}\|\vec{u}\|_{L^2_v
L^4_h(D)^3} \le C \|\vec{u}\|_{2}^{\frac{1}{2}} \| \nabla_h \vec{u}\|_{2}^{\frac{1}{2}}
\end{align}
\textbf{Proof.} See in \cite{paicuraugel}.
\end{Lem}


\begin{Lem}
Let   $g$  be a $2\pi$-periodic and mean zero   function with $g(x) \in H^s{([-\pi, \pi])}$, and let $s > \frac{1}{2} $, then we have the following inequality 
\begin{align}
\label{hanialiagmon}
\|g\|_{L^\infty} \le C\| g \|_{L^2}^{1-\frac{1}{2s}} \| g\|_{H^s}^{\frac{1}{2s}}.
\end{align}
\end{Lem}
\textbf{Proof.} The proof is given in Appendix.
\begin{Lem}
\label{hali}
Let $s > \frac{1}{2}$, and suppose that $\vec{u}$ and $ \partial^s \vec{u} \in L^2(D)^3$, then there exists a constant $C>0$
such that 
$$\|\vec{u}\|_{L^\infty _v
L^2_h(D)^3} \le C \| \vec{u} \|_{2}^{1-\frac{1}{2s}}\| \partial_3^{s} \vec{u} \|_{2}^{\frac{1}{2s}}. $$
\end{Lem}
\textbf{Proof.} We first apply inequality (\ref{hanialiagmon})   in the vertical
variable and then use the H\"{o}lder inequality in the horizontal variable.
We get,
\begin{equation}
\begin{split}
\sup_{x_3 \in [ -\pi, \pi]} \left( \int_{\Omega} |\vec{u}(\vec{x})|^2 dx_1 dx_2 \right)^{\frac{1}{2}}  \le \left( \int_{\Omega} \sup_{x_3 \in [ -\pi, \pi]} |\vec{u}(\vec{x})|^2 dx_1 dx_2 \right)^{\frac{1}{2}} \\
\quad \le  C \left( \int_{\Omega} \left(\int_{-\pi}^{\pi}  |\vec{u}(\vec{x})|^2 dx_3 \right)^{1-\frac{1}{2s}} \left( \int_{-\pi}^{\pi}  |\partial_3^s \vec{u}(\vec{x})|^2 dx_3\right)^{\frac{1}{2s}} dx_1 dx_2 \right)^{\frac{1}{2}}\\
\le  C\| \nabla\vec{u} \|_{2}^{1-\frac{1}{2s}} \| \partial_3^{s} \nabla\vec{u}\|_{2}^{\frac{1}{2s}}.
\end{split}
\end{equation}
The following lemma will be useful for the
estimate of the nonlinear term.
\begin{Lem}
\label{nonlinear}
There exists a positive constant $C >0$ such that, for any $s>\frac{1}{2}$ and for any  smooth
enough divergence-free vector fields $\vec{u},$   $\vec{v},$ and $\vec{w}$,
the following estimates hold,
$$ (i) \quad \left|\left( (\vec{u} \cdot\nabla) \vec{v}, \vec{w} \right)\right| \le C\|\vec{u}\|_{2}^{\frac{1}{2}} \| \nabla \vec{u}\|_{2}^{\frac{1}{2}}
 \|\nabla \vec{v} \|_{2}^{1-\frac{1}{2s}} \| \partial_3^{s} \nabla\vec{v}\|_{2}^{\frac{1}{2s}}  \|\vec{w}\|_{2}^{\frac{1}{2}} \| \nabla \vec{w}\|_{2}^{\frac{1}{2}},  $$
 and 
 $$ (ii) \quad \left|\left( (\vec{u} \cdot\nabla) \vec{v}, \vec{w} \right)\right| \le C\|\vec{u}\|_{2}^{\frac{1}{2}} \| \nabla \vec{u}\|_{2}^{\frac{1}{2}} \|\vec{v}\|_{2}^{\frac{1}{2}} \| \nabla \vec{v}\|_{2}^{\frac{1}{2}}
 \|\nabla \vec{w} \|_{2}^{1-\frac{1}{2s}} \| \partial_3^{s} \nabla\vec{w}\|_{2}^{\frac{1}{2s}}.  $$
\end{Lem}
\textbf{Proof.} 
(i) By using lemma \ref{paicuraugel} and lemma \ref{hali} we get 
 \begin{equation*}
\begin{split}
\left|\left( (\vec{u} \cdot\nabla) \vec{v}, \vec{w} \right)\right| \le  \|\vec{u}\|_{L^2_vL^4_h(D)^3} \|\nabla \vec{v}\|_{L^\infty _v
L^2_h(D)^3}  \|\vec{w}\|_{L^2_vL^4_h(D)^3}\\
 \le 
C \|\vec{u}\|_{2}^{\frac{1}{2}} \| \nabla \vec{u}\|_{2}^{\frac{1}{2}}
 \|\nabla\vec{v} \|_{2}^{1-\frac{1}{2s}} \| \partial_3^{s}\nabla\vec{v}\|_{2}^{\frac{1}{2s}} \|\vec{w}\|_{2}^{\frac{1}{2}} \| \nabla \vec{w}\|_{2}^{\frac{1}{2}}  
 \end{split}
 \end{equation*}
 (ii) We have  $$\left|\left( (\vec{u} \cdot\nabla) \vec{v}, \vec{w} \right)\right| \le  \displaystyle \left|  \int_{D} \vec{v} \otimes \vec{u} : \nabla \vec{w} \ d\vec{x} \right|, $$ thus,  in order to prove  the second inequality it is sufficient to swap  the roles of $\vec{v}$ and $\vec{w}$.
 
 \subsection{ The vertical filter and  the vertical deconvolution operator}
 Let  $\bv$  be a smooth function  of the form $\bv=\sum_{ {k_3}\in \Z\setminus\{0\}} \bc_{k_3}(x_1,x_2) \expk $.
The action of the vertical filter on $\bv(\vec{x})=\displaystyle \sum_{ {k_3}\in \Z\setminus\{0\}} \bc_{k_3}(x_1,x_2) \expk $ can be written as $ {\mathbb{A}_{3,\theta}({\bv})}=\displaystyle \sum_{ {k_3}\in \Z\setminus\{0\}} \mathcal{A}_{\theta}(k_3)\bc_{k_3}(x_1,x_2) \expk$, where the Fourier transform  with respect to $x_3 $ of the vertical filter is given by 
\begin{equation}
\mathcal{A}_{\theta}(k_3)=\left(1+\alpha^{2\theta}|{k_3}|^{2\theta}\right).
\end{equation}
Therefore, by using the Parseval's identity with respect to $x_3$ we get,
\BEQ
\| \mathbb{A}_{3,\theta}^{\frac{1}{2}} \bv \|_{2}^2= \| \bv  \|_{2}^2 + \alpha^{2\theta}\| \partial_3^{\theta}\bv  \|_{2}^2=\left( \mathbb{A}_{3,\theta} \bv, \bv \right).
\EEQ

Next, we recall some properties  of the vertical filter. These properties are proved in \cite{Berselli2012} for the horizontal filter and hold true for the vertical filter.   
\begin{Lem}
\label{luigi}
Let f be a smooth function and  $\vec{w}$  a smooth, space-periodic with respect to $x_3$, and divergence-free vector field defined on the domain D such that $w \cdot n  = 0$ on $\partial D$. Let the filtering be defined by (\ref{vertical}). 
 Then, the following identities  hold true:
 \begin{align}
\left( \overline{f}, \vec{w} \right)= \left( f, \overline{\vec{w}} \right),\label{luigi1} \\
\overline{\vec{w} - \alpha^{2\theta} \partial_3^{2\theta} \vec{w}} = \overline{\vec{w}} - \alpha^{2\theta} \partial_3^{2\theta}\overline{\vec{w}}, \label{luigi2}
\end{align}
and 
\begin{align} \left( \overline{\nabla \cdot(\vec{w} \otimes \vec{w} )}, \mathbb{A}_{3,\theta} \vec{w} \right)= \left(\nabla \cdot(\vec{w} \otimes \vec{w}) , \vec{w}\right) =0. \label{luigi3}
\end{align}
\end{Lem}
\textbf{Proof.} See in \cite{Berselli2012}. 

The deconvolution operator  $ D_{N,\theta} = \sum_{i = 0}^N ( I - \mathbb{A}_{3,\theta}^{-1} )^i$  is  constructed by using  the vertical filter  with fractional regularization  (\ref{vertical}). For a fixed $N >0$ and for $\theta =1$, we recover a vertical operator form from the Van Cittert deconvolution operator (see    \cite{AdamStolz} and   \cite{dunca06}).  
A straightforward calculation yields
\BEQ \begin{split} D_{N,\theta}  \left ( \sum_{ {k_3}\in \Z\setminus\{0\}} \bc_{k_3}(x_1,x_2) \expk \right ) \hskip 7cm \\
=  
\sum_{ {k_3}\in \Z\setminus\{0\}}\left(1+ \alpha^{2 \theta} | {k_3} |^{2 \theta}\right) \left ( 1- \left ( \frac{\alpha^{2 \theta} | {k_3} |^{2 \theta}}
 {1 +  \alpha^{2 \theta} | {k_3} |^{2 \theta}}\right )^{N+1} \right )  \bc_{k_3}(x_1,x_2) \expk.
\end{split} \EEQ 
Thus 
 \BEQ D_{N,\theta}  \left ( \sum_{ {k_3}\in \Z\setminus\{0\}} \bc_{k_3}(x_1,x_2) \expk \right ) =  \sum_{ {\bk}\in \Z\setminus\{0\}}  \mathcal{D}_{N,\theta}(k_3)  \bc_{k_3}(x_1,x_2) \expk, \EEQ 
 where we have for  ${k_3}\in \Z\setminus\{0\},$ and $ \theta \ge 0$,
 \begin{align}
 \mathcal{D}_{0,\theta}(k_3) &= 1, \\
 1 \le  \mathcal{D}_{N,\theta}(k_3) &\le N+1 \quad \hbox { for each } N >0, \\
 \hbox{ and  }   \mathcal{D}_{N,\theta}(k_3) &\le {\mathcal{A}_{3,\theta}}  \textrm{ for  a fixed } \alpha >0. \label{haniali}
  \end{align}
%
%
From the previous hypothesis, one can prove the following Lemma  by adapting the results summarized in the isotropic case
in \cite{bresslilewandowski}: 
\begin{Lem}
\label{lemme1}
For  $\theta \ge 0$, ${k_3}\in \Z\setminus\{0\}$ and for each $N >0$, there exists a constant $C >0$ such that for all $\bv$ sufficiently smooth  we have 
\begin{align}
  &\|\bv\|_{2} \le    \|D_{N,\theta} \left ( \bv \right )\|_{2} \le  (N+1)  \|\bv\|_{2},\label{lemme13a}\\
    &\|\bv\|_{2} \le   C \|D_{N,\theta} \left ( \bv \right )\|_{2} \le  C \|{\mathbb{A}_{3,\theta}}^{\frac{1}{2}}D_{N,\theta}^{\frac{1}{2}}(\bv)\|_{2},\label{lemme13b}\\
     &\|{\mathbb{A}_{3,\theta}}^{\frac{1}{2}}D_{N,\theta}^{\frac{1}{2}}(\overline{\bv})\|_{2}  \le  \|\bv\|_{2}, \label{lemme13}\\    
      &\|\bv\|_{2}^2 + \alpha^{2\theta}\| \partial_3^{\theta}\bv  \|_{2}^2 \le    \|{\mathbb{A}_{3,\theta}}^{\frac{1}{2}}D_{N,\theta}^{\frac{1}{2}}(\bv)\|_{2}^2.\label{lemme13d}
  \end{align}
\end{Lem}

\section{Existence and uniqueness results}

In this section, we give a definition of what is called a regular weak solution to problem (\ref{dalpha ns}).  Then, we give the proof of theorem \ref{TH1}.

\begin{Def}
Let  ${\bef} \in L^{2}(0,T;\vec{H})$ be a divergence free function and  $\bv_0 \in \vec{H}$.  For any $0 \le \theta \le 1$ 	and $ 0 \le N <\infty,$  $\bw$ is called a ``regular" weak solution to  problem (\ref{dalpha ns}) if the following
properties are satisfied:
\begin{align}
\bw &\in \mathcal{C}_{}(0,T;\vec{H}) \cap
L^2(0,T;\vec{V}),\label{bv12}\\
\partial_3^{\theta} \bw &\in \mathcal{C}_{}(0,T;\vec{H}) \cap
L^2(0,T;\vec{V}),\label{bv12bis}\\
\partial_t\bw&\in   L^{2}(0,T;\vec{V}^{*}),
\label{bvt}
\end{align}
and the velocity $\bw$ fulfills
\begin{equation}
\begin{split}
\int_0^{T}  \left( \partial_t \bw,   \bfi \right)  + \left(  \overline{ \nabla \cdot (D_{N,\theta}(\bw) \otimes  D_{N,\theta}(\bw)}), 
\bfi \right)  +  \nu \left( 
\nabla \bw, \nabla \bfi  \right)  \; dt\\
 =\int_0^{T}   \left( \overline{\bef}, \bfi \right)\; dt  
\qquad \textrm{ for all } {\bfi}\in L^{2}(0,T; \vec{V}).
\end{split}\label{weak1}
\end{equation}
Moreover,
\begin{equation}
\bw(0)= \bw_0.
\label{intiale}
\end{equation}
\end{Def}

\subsection{Proof of Theorem \ref{TH1}}
The proof of Theorem \ref{TH1} follows the classical scheme. 
 We start by constructing
approximated solution $\bw^n$ via Galerkin method. Then, we seek for a priori
estimates that are uniform with respect to $n$. Next, we take  the limit in the
equations after having used compactness properties. Finally, we show that the constructed solution is unique thanks to Gronwall's lemma \cite{RT83}.\\

\textbf{Step 1}(Galerkin approximation).
 We denote by 
$\left\{ \bfi^{k_3} \right\}_{|k_3|=1}^{+\infty }$ the eigenfunctions of the Stokes operator on D, with Dirichlet boundary conditions on $\partial D$ and with periodicity only  with respect to $x_3$. The explicit expression of these eigenfunctions can be found in 
\cite{Rummler1}.  
These eigenfunctions
are linear
combinations of $\mathcal{W}_{k_3}(x_1,x_2)\expk$, where $k_3 \in \Z\setminus\{0\},$ for certain  families of smooth functions $\mathcal{W}_{k_3} : \Omega \rightarrow \R$ vanishing at $\partial \Omega$.

We set 
\begin{equation}
\begin{split}
\bw^n(t,\bx)=\sum_{k_3 \in \Z\setminus\{0\}, |k_3| \le n }\bc_{k_3}^n (x_1, x_2, t) \bfi^{k_3}(\bx).
\end{split}
\end{equation}
 Since we are dealing with real
 functions, one has to take suitable conjugation and  alternatively use linear combinations of sines and cosines in the
 variable $x_3$.
We look for $\bw^n(t,\bx) $ which is determined by the following system of equations


\begin{equation}
\left\{
\begin{split}
 &\left( \partial_t\bw^n, \bfi^{k_3} \right)  + (\overline{\nabla \cdot( D_{N,\theta}(\bw^n)  \otimes D_{N,\theta}(\bw^n)) }, 
\bfi^{k_3}) +  \nu(
\nabla \bw^n, \nabla \bfi^{k_3}  )\; \\
& \qquad = \langle \bef, \bfi^{k_3} \rangle \; ,
\quad {|k_3|=1},2,...,n,\\
&{\bw^n(0)=P_n(\bw_0)=P_n(\overline{\bv_0})},
\end{split}\label{weak1galerkine}
\right.
\end{equation}
 where  $P_n$ denotes the projection operator over $H_n := \textrm{ Span}\langle\bfi^{1}, ... ,   \bfi^{n}\rangle$.
Therefore, the classical Caratheodory theory \cite{Wa70}   implies the short-time existence of solutions 
to (\ref{weak1galerkine}).  Next, we derive  estimates on $\bc^n$ that is uniform w.r.t. $n$.
These estimates  imply that the  solution of  (\ref{weak1galerkine}), constructed on a short time interval $[0, T^n[,$ exists for all $t \in [0, T]$.\\

\textbf{Step 2} (A priori estimates)
We need to derive an energy inequality for $\bw^n $. This can be obtained  by using $ \mathbb{A}_{3,\theta}D_{N,\theta}({\bw}^n)$ as a test function  in (\ref{weak1galerkine}). Thanks to the explicit expression  of the eigenfunctions,  and the properties  of $ \mathbb{A}_{3,\theta} $ and $D_{N,\theta}$, the quantity $ \mathbb{A}_{3,\theta}D_{N,\theta}({\bw}^n)$ is a legitimate test function since it still belongs to $H_n.$ Standard manipulations and the use of Lemma \ref{luigi}  combined with  the following identities
\begin{equation}
\begin{split}
 \label{divergencfreebar1}
\left(\partial_t{\bw}^n, \mathbb{A}_{3,\theta}D_{N,\theta}({\bw}^n)\right)=\frac{1}{2}\frac{d}{dt}\left(\|D_{N,\theta}^{\frac{1}{2}}({\bw}^n) \|_{2}^2 +\alpha^{2\theta}\|\partial_3^{\theta}D_{N,\theta}^{\frac{1}{2}}({\bw}^n) \|_{2}^2 \right)\\
 = \frac{1}{2}\frac{d}{dt}\|\mathbb{A}_{3,\theta}^{\frac{1}{2}}D_{N,\theta}^{\frac{1}{2}}({\bw}^n) \|_{2}^2,
\end{split}
\end{equation}
\begin{equation}
\begin{split}
 \label{divergencfreebar2}
\left(-\Delta{\bw}^n,   \mathbb{A}_{3,\theta}D_{N,\theta}^{}({\bw}^n)\right)=\left(\|\nabla D_{N,\theta}^{\frac{1}{2}}({\bw}^n)\|_{2}^2+ \|\partial_3^{\theta}\nabla D_{N,\theta}^{\frac{1}{2}}({\bw}^n)\|_{2}^2 \right) \\
 =
\|\nabla \mathbb{A}_{3,\theta}^{\frac{1}{2}}D_{N,\theta}^{\frac{1}{2}}({\bw}^n)\|_{2}^2,
\end{split}
\end{equation}
and 
\begin{equation}
 \label{divergencfreebar3}
\left( \overline{\bef}, \mathbb{A}_{3,\theta}D_{N,\theta}^{}({\bw}^n) \right)= \left(D_{N,\theta}^{\frac{1}{2}}({\bef}),D_{N,\theta}^{\frac{1}{2}}({\bw}^n) \right),
\end{equation}
 lead to the a priori estimate
\begin{equation}
\begin{array}{lllll}
 \label{apriori1}
\displaystyle \frac{1}{2}\| \mathbb{A}_{3,\theta}^{\frac{1}{2}}D_{N,\theta}^{\frac{1}{2}}({\bw}^n) \|_{2}^2 
+ \displaystyle \nu\int_{0}^{t}\| \nabla \mathbb{A}_{3,\theta}^{\frac{1}{2}}D_{N,\theta}^{\frac{1}{2}}({\bw}^n) \|_{2}^2  \ ds\\
 \quad = \displaystyle \int_{0}^{t} \left(D_{N,\theta}^{\frac{1}{2}}({\bef}), D_{N,\theta}^{\frac{1}{2}}({\bw}^n) \right) \ ds 
 +\displaystyle \frac{1}{2}\| \mathbb{A}_{3,\theta}^{\frac{1}{2}}D_{N,\theta}^{\frac{1}{2}}(\overline{\bv}^{n}_{0}) \|_{2}^2.
\end{array}
\end{equation}
By  using the Cauchy-Schwartz inequality, the Poincar\'{e} inequality combined with the Young inequality and inequality (\ref{lemme13}),   we conclude from  (\ref{apriori1}) the following inequality
 \begin{equation}
 \label{iciapriori12bis}
 \begin{split}
\sup_{t \in [0,T^n[}\| \mathbb{A}_{3,\theta}^{\frac{1}{2}}D_{N,\theta}^{\frac{1}{2}}({\bw}^n) \|_{2}^2 + \nu\int_{0}^{t} \| \nabla \mathbb{A}_{3,\theta}^{\frac{1}{2}}D_{N,\theta}^{\frac{1}{2}}({\bw}^n) \|_{2}^2  \ ds  \\
\le \| {\bv}^{n}_{0} \|_{2}^2 + \frac{C(N+1)}{\nu}\int_{0}^{T} \| {\bef} \|_{2}^2  \ ds
\end{split}
\end{equation}
which immediately implies that the existence time is independent of $n$ and it is possible to take $T=T^n$.\\ 
We deduce from (\ref{iciapriori12bis}) and (\ref{lemme13d}) that 
 \begin{equation}
 \label{iciapriori12}
 \begin{split}
\sup_{t \in [0,T]}\left(\|{\bw}^n\|_{2}^2+ \|\partial_3^{\theta}{\bw}^n\|_{2}^2 \right) + \nu\int_{0}^{t}\left(\| \nabla {\bw}^n\|_{2}^2+ \|\partial_3^{\theta}\nabla {\bw}^n\|_{2}^2 \right)  \ ds \\
 \le  \| {\bv}^{}_{0} \|_{2}^2 + \frac{C(N+1)}{\nu}\int_{0}^{T} \| {\bef} \|_{2}^2  \ ds.
\end{split}
\end{equation}
It follows from the above inequality that 
 \begin{align}
\label{vbar1sansdn}
 {\bw}^n \in L^{\infty}(0,T ; \vec{H}) \cap L^{2}(0,T ; \vec{V}),\\
  \partial^{\theta}_{3}{\bw}^n \in L^{\infty}(0,T ; \vec{H}) \cap L^{2}(0,T ; \vec{V}).
  \label{bisvbar1sansdnbis}
 \end{align}
 Thus, in one hand we get, 
  \begin{align}
\label{vbar1sansdnbis}
 \Delta {\bw}^n \in  L^{2}(0,T ; \vec{V}^{*}),\\
  \partial^{\theta}_{3} \Delta {\bw}^n \in  L^{2}(0,T ; \vec{V}^{*}). \label{vbar1sansdnbisbis}
 \end{align}
 In the other hand, it follows from (\ref{lemme13a}) that 
 \begin{align}
\label{icivbar1}
 D_{N,\theta}({\bw}^n) \in L^{\infty}(0,T ; \vec{H}) \cap L^{2}(0,T ; \vec{V}),\\
  \partial^{\theta}_{3} D_{N,\theta}({\bw}^n) \in L^{\infty}(0,T ; \vec{H}) \cap L^{2}(0,T ; \vec{V}).\label{icivbar1bis}
 \end{align}
 The inequality (\ref{icivbar1}) allows us to find estimates on the nonlinear term $\overline{\nabla \cdot  (D_{N,\theta}({\bw}^n) \otimes  D_{N,\theta}({\bw}^n))}.$  For that 
 we use lemma \ref{nonlinear}. 
For all  $\bfi \in  \vec{V}$ we have 
\begin{equation}
\begin{split}
\left| \left(   \overline{\nabla \cdot  (D_{N,\theta}({\bw}^n) \otimes  D_{N,\theta}({\bw}^n))}, \bfi \right) \right|  \le \left| \left( {\nabla \cdot  (D_{N,\theta}({\bw}^n) \otimes  D_{N,\theta}({\bw}^n))}  , \overline{\bfi} \right) \right| \\ 
 \le C\|D_{N,\theta}({\bw}^n)\|_2\| \nabla D_{N,\theta}({\bw}^n)\|_2
 \|\nabla \overline{\bfi}\|_{2}^{1-\frac{1}{2\theta}} \| \nabla  \partial_3^{\theta} \overline{\bfi} \|_{2}^{\frac{1}{2\theta}}. 
 \end{split}
\end{equation}
We also have 
\begin{equation}
\begin{split}
\left| \left(  \mathbb{A}_{3,\theta}^{\frac{1}{2}} \overline{\nabla \cdot  (D_{N,\theta}({\bw}^n) \otimes  D_{N,\theta}({\bw}^n))}, \bfi \right) \right|  \le \left| \left( {\nabla \cdot  (D_{N,\theta}({\bw}^n) \otimes  D_{N,\theta}({\bw}^n))}  , \mathbb{A}_{3,\theta}^{\frac{1}{2}}\overline{\bfi} \right) \right| \\ 
 \le C\|D_{N,\theta}({\bw}^n)\|_2\| \nabla D_{N,\theta}({\bw}^n)\|_2
 \|\nabla\mathbb{A}_{3,\theta}^{\frac{1}{2}} \overline{\bfi}\|_{2}^{1-\frac{1}{2\theta}} \| \nabla  \partial_3^{\theta}\mathbb{A}_{3,\theta}^{\frac{1}{2}} \overline{\bfi} \|_{2}^{\frac{1}{2\theta}} 
 \end{split}
\end{equation}
Since $ \| \nabla  \partial_3^{\theta} \overline{\bfi} \|_{2}^{\frac{1}{2\theta}} \le  \| \nabla {\bfi} \|_{2}^{\frac{1}{2\theta}},$  $ \|\nabla\mathbb{A}_{3,\theta}^{\frac{1}{2}} \overline{\bfi}\|_{2}^{1-\frac{1}{2\theta}}  \le \|\nabla{\bfi}\|_{2}^{1-\frac{1}{2\theta}}$
and $\| \nabla  \partial_3^{\theta}\mathbb{A}_{3,\theta}^{\frac{1}{2}} \overline{\bfi} \|_{2}^{\frac{1}{2\theta}}  \le \| \nabla {\bfi} \|_{2}^{\frac{1}{2\theta}} $ we get that 

\begin{equation}
\begin{split}
\left| \left(  \overline{\nabla \cdot  (D_{N,\theta}({\bw}^n) \otimes  D_{N,\theta}({\bw}^n))}, \bfi \right) \right|  
 \le C\|D_{N,\theta}({\bw}^n)\|_2\| \nabla D_{N,\theta}({\bw}^n)\|_2 \|\nabla{\bfi}\|_{2},
 \end{split}
\end{equation}
and 
\begin{equation}
\begin{split}
\left| \left(  \mathbb{A}_{3,\theta}^{\frac{1}{2}} \overline{\nabla \cdot  (D_{N,\theta}({\bw}^n) \otimes  D_{N,\theta}({\bw}^n))}, \bfi \right) \right|  
 \le C\|D_{N,\theta}({\bw}^n)\|_2 \| \nabla D_{N,\theta}({\bw}^n)\|_2 \|\nabla{\bfi}\|_{2}.
 \end{split}
\end{equation}

Thus we obtain 
\begin{align}
\overline{\nabla \cdot  (D_{N,\theta}({\bw}^n) \otimes  D_{N,\theta}({\bw}^n))} \in  L^{2}(0,T ; \vec{V}^{*}),\\
\mathbb{A}_{3,\theta}^{\frac{1}{2}} \overline{\nabla \cdot  (D_{N,\theta}({\bw}^n) \otimes  D_{N,\theta}({\bw}^n))} \in  L^{2}(0,T ; \vec{V}^{*}).\label{ici17}
\end{align}


From eqs. (\ref{weak1galerkine}), (\ref{vbar1sansdn}), (\ref{ici17}), we also obtain that 
 \begin{equation}
\label{icivtemps}
 \partial_t \bw^n \in L^{2}(0,T ; \vec{V}^{*}), 
\end{equation}
and 
 \begin{equation}
\label{icivtempsbis}
 \partial_t \mathbb{A}_{3,\theta}^{\frac{1}{2}} \bw^n \in L^{2}(0,T ; \vec{V}^{*}). 
\end{equation}

\textbf{Step 3} (Limit $n \rightarrow \infty$) It follows from the estimates (\ref{vbar1sansdn})-(\ref{icivtempsbis}) and the Aubin-Lions compactness lemma
(see \cite{sim87} for example) that there exists  a  not relabeled  subsequence of $\bw^n$ and $\bw$ such that
\begin{align}
\bw^n &\rightharpoonup^* \bw &&\textrm{weakly$^*$ in } L^{\infty}
(0,T;\vec{H}), \label{c122}\\
D_{N,\theta}(\bw^n) &\rightharpoonup^* D_{N,\theta}(\bw) &&\textrm{weakly$^*$ in } L^{\infty}
(0,T;\vec{H}), \label{DNc122}\\
\bw^n &\rightharpoonup \bw &&\textrm{weakly in }
L^2(0,T;\vec{V}), \label{nc22}\\
\partial_{3}^{\theta}\bw^n &\rightharpoonup \partial_{3}^{\theta}\bw &&\textrm{weakly in }
L^2(0,T;\vec{V}), \label{partialNc22}\\
D_{N,\theta}(\bw^n) &\rightharpoonup D_{N,\theta}(\bw) &&\textrm{weakly in }
L^2(0,T;\vec{V}), \label{DNnc22}\\
\partial_{3}^{\theta}D_{N,\theta}(\bw^n) &\rightharpoonup \partial_{3}^{\theta}D_{N,\theta}(\bw) &&\textrm{weakly in }
L^2(0,T;\vec{V}), \label{partialDNnc22}\\
\partial_t\bw^n&\rightharpoonup \partial_t\bw &&\textrm{weakly in } L^{2}
(0,T;\vec{V}^{*}),
\label{nc322}\\
 \partial_t \mathbb{A}_{3,\theta}^{\frac{1}{2}} \bw^n &\rightharpoonup \partial_t \mathbb{A}_{3,\theta}^{\frac{1}{2}}\bw &&\textrm{weakly in } L^{2}(0,T;\vec{V}^{*}),\\
\bw^n &\rightarrow \bw &&\textrm{strongly in  }
L^2(0,T;\vec{H}),
\label{c83icil2}\\
D_{N,\theta}(\bw^n) &\rightarrow D_{N,\theta}(\bw) &&\textrm{strongly in  }
L^2(0,T;\vec{H}),
\label{DNc83icil2}
\end{align}
 From (\ref{DNnc22}) and (\ref{DNc83icil2}),  it follows   that  
 
 \begin{align}
 \overline{\nabla\cdot(D_{N,\theta}({\bw}^n)  \otimes D_{N,\theta}({\bw}^n))} &\rightarrow \overline{\nabla \cdot(D_{N,\theta}({\bw}) \otimes  D_{N,\theta}(\bw))} &&\textrm{strongly in  }
L^1(0,T;L^{1}(D)^{3}),\label{c82int}
 \end{align}
Finally, since the
sequence $  \overline{\nabla \cdot (D_{N,\theta}({\bw}^n) \otimes D_{N,\theta}({\bw}^n))}$  is bounded in $L^2(0,T;\vec{V}^{*})$, it converges weakly, up to a
subsequence, to some $\chi $ in  $L^2(0,T;\vec{V}^{*})$. The previous result  and the  uniqueness of the limit
allow us to claim that $\chi=\overline{ \nabla \cdot(D_{N,\theta}(\bw) \otimes  D_{N,\theta}({\bw}))}$. Consequently,
 \begin{align}
     \overline{\nabla \cdot (D_{N,\theta}({\bw}^n) \otimes D_{N,\theta}({\bw}^n))} &\rightharpoonup  \overline{ \nabla \cdot(D_{N,\theta}(\bw) \otimes  D_{N,\theta}({\bw}))}&&\textrm{weakly in }
 L^2(0,T;\vec{V}^{*}). \label{c22primen}\    
  \end{align}
The above established convergences are clearly sufficient to take the limit in (\ref{weak1galerkine})  and  conclude that  the velocity   $ \bw$  satisfies (\ref{weak1}). 
Moreover, 
from (\ref{nc22}) and (\ref{nc322}),
 one  can deduce by a classical interpolation  argument \cite{LiMa68}   that 
 \begin{equation}
 \bw \in  \mathcal{C}(0,T;\vec{H}).
\end{equation}
Furthermore, from  the  strong continuity of $\bw$ with respect to the time and with values in $\vec{H}$,  we deduce   that $\bw(0)=\bw_0$.\\

 \textbf{Step 4} (Uniqueness)
Next, we will show the continuous dependence of the  solutions on the initial data and in particular the uniqueness.\\
Let $ \theta > \frac{1}{2}$  and let $({\bw_1,q_1})$ and $({\bw_2,q_2})$   be any two solutions of (\ref{dalpha ns}) on the interval $[0,T]$, with initial values $\bw_1(0)$ and $\bw_2(0)$. Let us denote by  $\delta \vec{w}_{} =\bw_2-\bw_1$.
 We subtract the equations for $\bw_1$ from the equations for $\bw_2$. Then,
 we obtain
\begin{equation}
\label{matin1}
\begin{split}
\partial_{t} \delta\vec{w} -\nu \Delta \delta\vec{w} +   \overline{\nabla\cdot(D_{N,\theta}(\bw_2) \otimes D_{N,\theta}(\bw_2)}) - \overline{\nabla \cdot (D_{N,\theta}(\bw_1) \otimes D_{N,\theta}(\bw_1})) =0,
\end{split}
\end{equation}
and $\delta\vec{w} = 0$  at the initial time.\\ 
Applying $\mathbb{A}_{3,\theta}^{\frac{1}{2}}  $  to  (\ref{matin1})  we obtain

\begin{equation}
\label{matin1ha}
\begin{split}
\mathbb{A}_{3,\theta}^{\frac{1}{2}}  \partial_{t} \delta\vec{w} -\nu \mathbb{A}_{3,\theta}^{\frac{1}{2}} \Delta \delta\vec{w} +  \mathbb{A}_{3,\theta}^{\frac{1}{2}} \overline{\nabla\cdot(D_{N,\theta}(\bw_2) \otimes D_{N,\theta}(\bw_2)})\\
- \mathbb{A}_{3,\theta}^{\frac{1}{2}} \overline{\nabla (\cdot D_{N,\theta}(\bw_1) \otimes D_{N,\theta}(\bw_1}))=0
\end{split}
\end{equation}
One can take $ \mathbb{A}_{3,\theta}^{\frac{1}{2}} D_{N,\theta}(\delta\vec{w}) \in {L}^2(0,T;\vec{V})$
 as test function  in   (\ref{matin1ha}). Let us mention that, $ \partial_{t} \mathbb{A}_{3,\theta}^{\frac{1}{2}} \delta\vec{w} \in L^{2}(0,T ; \vec{V}^{*}) $  and $\mathbb{A}_{3,\theta}^{\frac{1}{2}} D_{N,\theta}(\delta\vec{w}) \in {L}^2(0,T;\vec{V})$. Thus,   by using Lions-Magenes Lemma \cite{LiMa68}  we have
 $$
 \langle  \partial_{t} \mathbb{A}_{3,\theta}^{\frac{1}{2}} \delta\vec{w},  \mathbb{A}_{3,\theta}^{\frac{1}{2}} D_{N,\theta}(\delta\vec{w}) \rangle_{\vec{V}^{*}, \vec{V} }= \frac{1}{2}\frac{d}{d}\| \mathbb{A}_{3,\theta}^{\frac{1}{2}}D_{N,\theta}^{\frac{1}{2}}(\delta\vec{w})\|_{2}^2.$$
Using lemma \ref{luigi} and the divergence free condition  we get:
 \begin{equation}
 \label{notuniformzerobis}
 \begin{split}
  \displaystyle 
 \frac{1}{2}\frac{d}{dt} \|\mathbb{A}_{3,\theta}^{\frac{1}{2}}D_{N,\theta}^{\frac{1}{2}}{(\delta \vec{w})}\|_{2}^{2} +\nu \|\nabla \mathbb{A}_{3,\theta}^{\frac{1}{2}}D_{N,\theta}^{\frac{1}{2}}{(\delta \vec{w})}_{}\|_{2}^{2} 
   \le \displaystyle \left(({D_{N,\theta}(\delta \vec{w})} \cdot \nabla) D_{N,\theta}({\bw}_{2}) ,  D_{N,\theta}^{} (\delta \vec{w}) \right).
   \end{split}
\end{equation}
We estimate now the right-hand side by using lemma \ref{nonlinear} as follows,
\begin{equation}
 \label{notuniformzero}
 \begin{split}
  \displaystyle 
 \displaystyle \left|\left(({D_{N,\theta}(\delta \vec{w})} \cdot \nabla)  D_{N,\theta}({\bw}_{2}) ,  D_{N,\theta}^{} (\delta \vec{w}) \right)\right| \hskip 6cm \\
\le \|D_{N,\theta}(\delta \vec{w})\|_{2}\|\nabla D_{N,\theta}(\delta \vec{w})\|_{2} \| \nabla  D_{N,\theta}(\vec{w}_{2})\|_{2}^{1-\frac{1}{2\theta}} \|\partial_{3}^{\theta} \nabla  D_{N,\theta}(\vec{w}_{2})\|_{2}^{\frac{1}{2\theta}}
   \end{split}
\end{equation}
Hence, by using  the Young inequality combined with lemma \ref{lemme1}, we obtain that there exists a constant $C > 0 $  that depends on $N$ and $\nu$ such 
that

\begin{equation}
 \label{notuniformzerohani}
 \begin{split}
  \displaystyle 
 \displaystyle \left|\left(({D_{N,\theta}(\delta \vec{w})} \nabla ) D_{N,\theta}({\bw}_{2}) ,  D_{N,\theta}^{} (\delta \vec{w}) \right)\right| \hskip 4cm\\
\le C \| \delta \vec{w}\|_{2}^2 \| \nabla \vec{w}_{2}\|_{2}^{2-\frac{1}{\theta}} \|\partial_{3}^{\theta} \nabla \vec{w}_{2}\|_{2}^{\frac{1}{\theta}} +\frac{\nu}{2} \| \nabla \delta \vec{w}\|_{2}^2
   \end{split}
\end{equation}
By using  (\ref{lemme13d}) we have 
\begin{equation}
\label{notuniform}
 \begin{split}
  \displaystyle 
   \frac{1}{2}\frac{d}{dt}\left( \|{\delta \vec{w}}\|_{2}^{2} +  \alpha^{2\theta} \|\partial_{3}^{\theta}{\delta \vec{w}}\|_{2}^{2} \right) + \nu \left( \|\nabla {\delta \vec{w}}_{}\|_{2}^{2} + \alpha^{2\theta} \|\nabla {\delta \vec{w}}_{}\|_{2}^{2} \right) \\
\le 
  \frac{1}{2}\frac{d}{dt} \|\mathbb{A}_{3,\theta}^{\frac{1}{2}}D_{N,\theta}^{\frac{1}{2}}({\delta \vec{w}})\|_{2}^{2} +\nu \|\nabla \mathbb{A}_{3,\theta}^{\frac{1}{2}}D_{N,\theta}^{\frac{1}{2}}{(\delta \vec{w})}_{}\|_{2}^{2} \end{split}
\end{equation}
From (\ref{notuniformzero}), (\ref{notuniform}) we get 
\begin{equation}
 \begin{split}
  \displaystyle 
 \displaystyle 
   \frac{d}{dt}\left( \|{\delta \vec{w}}\|_{2}^{2} +  \alpha^{2\theta} \|\partial_{3}^{\theta}{\delta \vec{w}}\|_{2}^{2} \right) + \nu \left( \|\nabla {\delta \vec{w}}_{}\|_{2}^{2} + \alpha^{2\theta} \|\nabla {\delta \vec{w}}_{}\|_{2}^{2} \right)   \\
 \le \displaystyle  C \| \delta \vec{w}\|_{2}^2 \| \nabla \vec{w}_{2}\|_{2}^{2-\frac{1}{\theta}} \|\partial_{3}^{\theta} \nabla \vec{w}_{2}\|_{2}^{\frac{1}{\theta}}.
\end{split}
\end{equation}
Since $ \| \nabla \vec{w}_{2}\|_{2}^{2-\frac{1}{\theta}} \|\partial_{3}^{\theta} \nabla \vec{w}_{2}\|_{2}^{\frac{1}{\theta}}  \in L^1([0,T])$, we conclude by  using Gronwall's inequality
the continuous dependence of the solutions on the initial data in the $L^{\infty}([0,T],\vec{H})$  norm. In particular, if ${\delta \vec{w}}^{}_{0}=0$ then ${\delta \vec{w}}=0$ and the solutions are unique for all $t \in [0,T] .$ In addition, since $T>0$ is arbitrary chosen, this solution may be uniquely extended for all time.\\ 
This finishes the proof of Theorem \ref{TH1}.\\

\textbf{Acknowledgement}
I would like to thank Taoufik Hmidi  for his valuable discussions, help and advices.\\

\textbf{Appendix:  Proof of the functional estimate (\ref{hanialiagmon})}

Let   $g$  be a $2\pi$-periodic and mean zero   function with $g(x) \in H^s{([-\pi, \pi])}$, and let $s > \frac{1}{2} $. If we  write $g(x)$ as the sum of its Fourier series $g(x)=\displaystyle \sum_{ {k}\in \Z\setminus\{0\}} g_k {e^{ {\rm i} \, { k} x}} $, then, we can estimate $\| g \|_{L^\infty }$ by 
\begin{equation*}
\| g \|_{L^\infty } \le \sum_{ {k}\in \Z\setminus\{0\}} |g_k|.
\end{equation*}
We then break up the sum into low and high wave-number components,
\begin{equation*}
\| g \|_{L^\infty } \le \sum_{ 0 <|k| \le \kappa } |g_k | +  \sum_{ |k| > \kappa } |g_k |.
\end{equation*}
We now use the Cauchy-Schwarz  inequality on each part,
\begin{equation*}
\begin{split}
\| g \|_{L^\infty } \le \left(\sum_{ 0 <|k| \le \kappa } |g_k |^2\right)^{\frac{1}{2}} \left(\sum_{ 0 <|k| \le \kappa } 1 \right)^{\frac{1}{2}}  + \left( \sum_{ |k| > \kappa } |k|^{2s}|g_k|^2 \right)^{\frac{1}{2}}\left( \sum_{ |k| > \kappa } |k|^{-2s}\right)^{\frac{1}{2}}. 
 \end{split}
\end{equation*}
Since 
\begin{equation*}
\begin{split}
\sum_{ 0 <|k| \le \kappa } 1  \le   C\kappa \qquad \textrm{ and } \qquad  \sum_{ |k| > \kappa } |k|^{-2s} \le  C\kappa^{-2s+1},
 \end{split}
\end{equation*}
the above inequality becomes 
\begin{equation*}
\begin{split}
\| g \|_{L^\infty } \le C \left( \kappa^{\frac{1}{2}} \|g\|_{L^2}  + \kappa^{-s+\frac{1}{2}}\|g\|_{H^s} \right).
 \end{split}
\end{equation*}
To make both terms on the right-hand side the same, we choose 
$$ \kappa =  \frac{\|g\|_{H^s}^{\frac{1}{s}}}{\|g\|_{L^2}^{\frac{1}{s}}}.$$
This yields the following estimate
\begin{align*}
\|g\|_{L^\infty} \le C\| g \|_{L^2}^{1-\frac{1}{2s}} \| g\|_{H^s}^{\frac{1}{2s}},
\end{align*}
which is (\ref{hanialiagmon}).


\end{document}